\documentclass[reqno, a4paper,12pt]{amsart}

\usepackage[capposition=bottom]{floatrow}
\usepackage{euscript,eufrak,verbatim}
\usepackage{graphicx}
\usepackage[usenames]{color}
\usepackage[colorlinks,linkcolor=red,anchorcolor=blue,citecolor=blue]{hyperref}
\usepackage{amsmath, mathtools}
\usepackage{mathtools}
\usepackage{amsthm}
\usepackage[all]{xy}
\usepackage{amssymb} 
\usepackage{bm}

\usepackage[all]{xy}

\usepackage{mathrsfs}
\usepackage{amscd}

\usepackage{enumitem}
\usepackage[numbers]{natbib}

\usepackage{mathptmx}
\usepackage[T1]{fontenc}

\usepackage{lipsum}

\makeatletter
\numberwithin{equation}{section}

\theoremstyle{cuplain}
\newtheorem{main theorem}{Main Theorem}
\newtheorem{theorem}{Theorem}[section]

\newtheorem{proposition}[theorem]{Proposition}
\newtheorem{claim}[theorem]{Claim}

\newtheorem*{theorem*}{Theorem}
\theoremstyle{definition}
\newtheorem{definition}[theorem]{Definition}
\newtheorem{remark}[theorem]{Remark}
\newtheorem{example}[theorem]{Example}

\newtheorem*{example*}{Example}
\newtheorem*{remark*}{Remark}

\newtheoremstyle{break}
  {\topsep}{\topsep}%
  {\itshape}{}%
  {\bfseries}{}%
  {\newline}{}%
\theoremstyle{break}

\newtheoremstyle{break}
  {\topsep}{\topsep}%
  {\normalshape}{}%
  {\bfseries}{}%
  {\newline}{}%
\newtheorem{breakdefinition}[theorem]{Definition}

\numberwithin{equation}{section}
\newcommand{\spa}{\hspace{1pt}}

\newcommand{\flo}{\mathscr}

\newcommand{\norm}[1]{\left\lVert#1\right\rVert}
\newcommand{\setcond}{\hspace{2pt} \middle| \hspace{2pt}}

\DeclareFontFamily{U}{stix2bb}{\skewchar\font127 }
\DeclareFontShape{U}{stix2bb}{m}{n} {<-> stix2-mathbb}{}
\DeclareMathAlphabet{\mathbb}{U}{stix2bb}{m}{n}

\begin{document}


\newcommand\titlelowercase[1]{\texorpdfstring{\lowercase{#1}}{#1}}

\font\mathptmx=cmr12 at 12pt


\title[\fontsize{13}{12}\mathptmx {\it{E\titlelowercase{xact} H\titlelowercase{ausdorff dimension of some sofic self-affine fractals}}}]{\Huge E\titlelowercase{xact} H\titlelowercase{ausdorff dimension} \protect{\\[9pt]} \titlelowercase{of some sofic self-affine fractals}}


\author[\fontsize{13}{12}\mathptmx {\it{N\titlelowercase{ima} A\titlelowercase{libabaei}}}]{\fontsize{13}{12}\mathptmx Nima Alibabaei}



\maketitle

\begin{abstract}
Previous work has shown that the Hausdorff dimension of sofic affine-invariant sets is expressed as a limit involving intricate matrix products. This limit has typically been regarded as incalculable. However, in several highly non-trivial cases, we demonstrate that the dimension can in fact be calculated explicitly. Specifically, the dimension is expressed as the solution to an infinite-degree equation with explicit coefficients, which also corresponds to the spectral radius of a certain linear operator. Our result provides the first non-trivial calculation of the exact Hausdorff dimension of sofic sets in $\mathbb{R}^3$. This is achieved by developing a new technique inspired by the work of Kenyon and Peres (1998).
\end{abstract}

\section{Introduction and results} \label{section: introduction}

\subsection{Overview} 
\begin{figure}[h!]
\includegraphics[width=15.8cm]{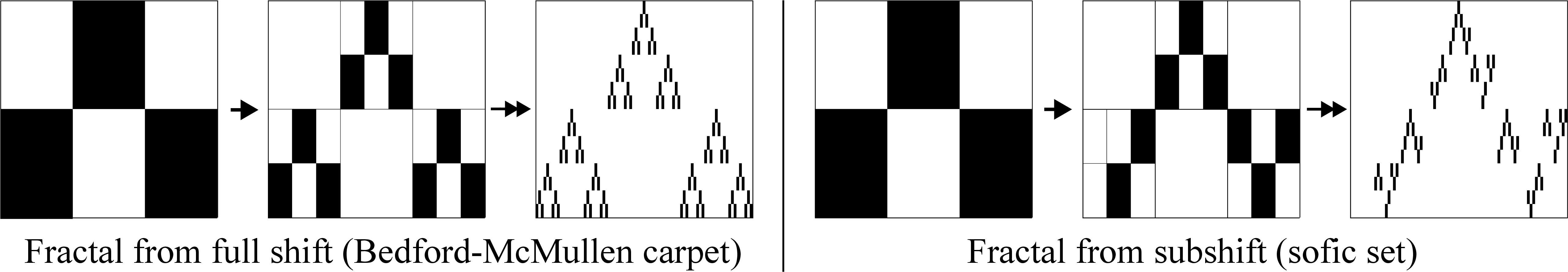}
\caption{The first few generations of self-affine fractals} \label{figure: self-affine fractals}
\end{figure}

\hfill\\
\quad Self-affine fractals, such as Bedford-McMullen carpets, have been extensively studied. Figure \ref{figure: self-affine fractals} illustrates the construction of these sets. Consider integers $1 < m_1 \leq m_2$, and partition the unit square into $m_1 m_2$ numbers of congruent rectangles (with $m_1 = 2$ and $m_2 = 3$ in Figure \ref{figure: self-affine fractals}). We then select some of these rectangles and repeat the process indefinitely, choosing which rectangles to retain at each step based on the initial selection. This process generalizes to higher Euclidean dimensions, and the resulting fractals are known as {\it{\textbf{self-affine sponges}}}. The Hausdorff dimension and the Minkowski (box) dimension of these sets are well understood \cite{Kenyon--Peres}. The underlying structure of these sets involves the symbolic dynamical systems known as full shifts, an infinite Cartesian product of finite symbols: $\{1, 2, \ldots, n\}^{\mathbb{N}}$. The shift-invariant nature of these systems ensures the affine-invariant (or "fractal") property of the self-affine sponges.

Next, we consider a subshift, a subset of a full shift that is shift-invariant. Using it, we can define a fractal in Euclidean space in a manner analogous to the construction of self-affine sponges. A particularly interesting class of subshifts is {\it{\textbf{sofic systems}}}, which arise from finite directed graphs. The fractals associated with sofic systems are referred to as {\it{\textbf{sofic sets}}} (see Figure \ref{figure: self-affine fractals}). 

 Sofic sets in $\mathbb{R}^2$ were studied in \cite{Kenyon--Peres: sofic}, where the Hausdorff dimension is expressed as a limit involving complicated matrix products. They introduced clever methods to calculate the Hausdorff dimension from this expression for various examples. Despite these insightful results, no further investigations have been made into the estimate of exact dimension of sofic sets.

This paper aims to derive an exact expression for the Hausdorff dimension of certain sofic sets in $\mathbb{R}^2$ and $\mathbb{R}^3$. We address the complexities that arise when the matrices lack the simplifying structures, such as commutativity or shared eigenvectors. By analyzing an example presented in \cite{Kenyon--Peres: sofic}, we introduce a technique we call {\it{tower-decomposition}} to break down the matrix product. For sofic sets in $\mathbb{R}^3$, this description requires operators defined on a space indexed by a tree, phenomena not encountered in planar cases. The Hausdorff dimension will be expressed as a solution to an equation of infinite degree, which coincides with the spectral radius of a certain linear operator. Our result provides the first exact calculation of the Hausdorff dimension for non-trivial sofic sets in $\mathbb{R}^3$, an achievement that has generally been considered highly challenging.

\subsection{Results in $\mathbb{R}^2$} \label{subsection: results in 2-dim} \hfill\\
\quad In \S \ref{subsection: results in 2-dim} and \S \ref{subsection: results in 3-dim}, we briefly outline our main results. The precise formulations of the main theorems require some preparation, so here we present only the ``rough versions'' of the statements along with illustrating examples. The rigorous statements and proofs of the theorems will be provided in \S \ref{section: theorems and proofs}, and the detailed calculations for the examples are given in \S \ref{section: examples}.

We now turn to the construction of sofic sets. Let $I$ be a set of labels. For instance, $I = \{0, 1\} \times \{0, 1, 2\}$ in Figure \ref{figure: self-affine fractals}. We begin by dividing the square into congruent rectangular pieces, each labeled with an element from $I$. These rectangles are then recursively subdivided into smaller rectangles, with the same labeling scheme applied at each stage. This process continues inductively, so that each element of $I^{\mathbb{N}}$ corresponds to a point within the square. It is important to note that this correspondence is not injective.

\begin{figure}[h!]
\includegraphics[width=15cm]{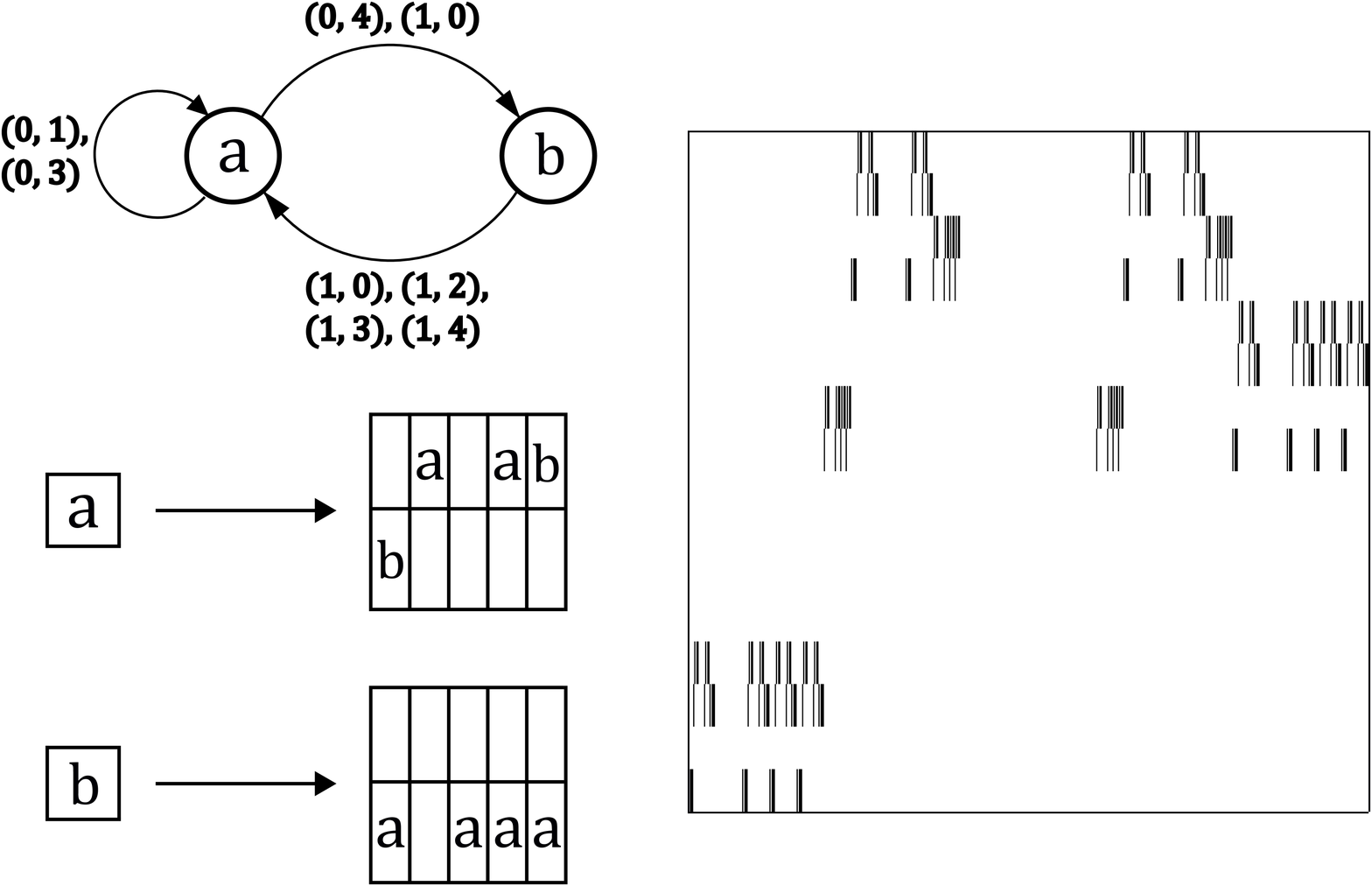}
\caption{A sofic set in $\mathbb{R}^2$ with $I = \{0,1\} \times \{0,1,2,3,4\}$} \label{figure: a sofic set}
\end{figure}

Next, we consider a directed graph $G$ where each edge is labeled with an element from $I$. An infinite path in $G$ corresponds to an element in $I^{\mathbb{N}}$, and we define $S \subset I^{\mathbb{N}}$ as the set of all such paths. The set $S$ is referred to as a sofic system. (We note that an additional condition on $G$ is required, which we omit here.) Via the aforementioned correspondence, $S$ defines a self-affine fractal in Euclidean space, commonly known as a sofic set. An example of such a set, with $I = \{0, 1\} \times \{0, 1, 2, 3, 4\}$, is illustrated in Figure \ref{figure: a sofic set}. (For a precise definition, see \S \ref{section: background}.)

Let $\mathbb{N}_0 = \{0, 1, 2, \ldots\}$ and $\mathrm{M}_n(\mathbb{N}_0)$ be the set of $n \times n$ matrices with components in $\mathbb{N}_0$. Consider integers $1 < m_1 \leq m_2$. Let $I = \{0, 1, \ldots, m_1-1\} \times \{0, 1, \ldots, m_2-1\}$ be the set of indices, and let $X \subset \mathbb{R}^2$ be the resulting sofic set. Let $I_1 = \{0, 1, \ldots, m_1-1\}$. By \cite[Theorem 3.2]{Kenyon--Peres: sofic}, there is an ``adjacency matrix'' $A_i \in \mathrm{M}_n(\mathbb{N}_0)$ for each $i \in I_1$ which yields the following formula of the Hausdorff dimension $\mathrm{dim_H}(X)$.
\begin{equation*}
\mathrm{dim_H}(X)
=
\lim_{N \to \infty} \spa \frac{1}{N} \log_{m_1^{}}{ \sum_{
				(u_1^{}, \ldots, u_N^{}) \in I_1^{N}}
			{\norm{A_{u_1^{}} \cdots A_{u_N^{}}}
			}^{\alpha}}.
\end{equation*}
Here, $\alpha = \log_{m_2}{m_1}$ and the norm of the matrices is arbitrary.

A matrix $A$ is said to be \textbf{primitive} if there is an integer $d$ such that every entry in $A^d$ is positive. Denote by $x^{\top}$ the transpose of a column vector $x \in \mathbb{R}^n$. We say that $A$ has a $1$-dimensional image when its image $\left\{ x^{\top} \spa A \setcond x \in \mathbb{R}^n \right\}$ is spanned by a single non-zero vector. (We consider matrices to be acting on $\mathbb{R}^n$ from right.) The following is the main result for sofic sets in $\mathbb{R}^2$. A rigorous statement of this result is given in Theorem \ref{theorem: tower decomposition for 2d}.

\begin{theorem*}
Suppose $\sum_{i \in I_1} A_i$ is primitive. Also assume that we have a natural number $L$ and a string $s = (s_1^{}, \ldots, s_L^{}) \in I_1^L$ such that the corresponding matrix $A_{s_1^{}} \cdots A_{s_L^{}}$ has a $1$-dimensional image. Then, there is an ``explicitly calculable'' constant $C_k \geq 0$ for each non-negative integer $k$ such that
\[\mathrm{dim}_{\mathrm{H}}(X) = \log_{m_1^{}}{r}, \]
where $r$ is the unique positive solution to the equation:
\[ r^L = C_0 + \frac{C_1}{r^1} + \frac{C_2}{r^2} + \cdots. \]
\end{theorem*}

This theorem provides an exact value for the Hausdorff dimension, as illustrated by the following examples.

\begin{example*}

Let $I = \{0, 1\} \times \{0, 1, 2\}$. Consider the directed graph in Figure \ref{figure: digraph for planar 01} labeled with $I$.
\begin{figure}[h!]
\includegraphics[width=\textwidth-5.5cm]{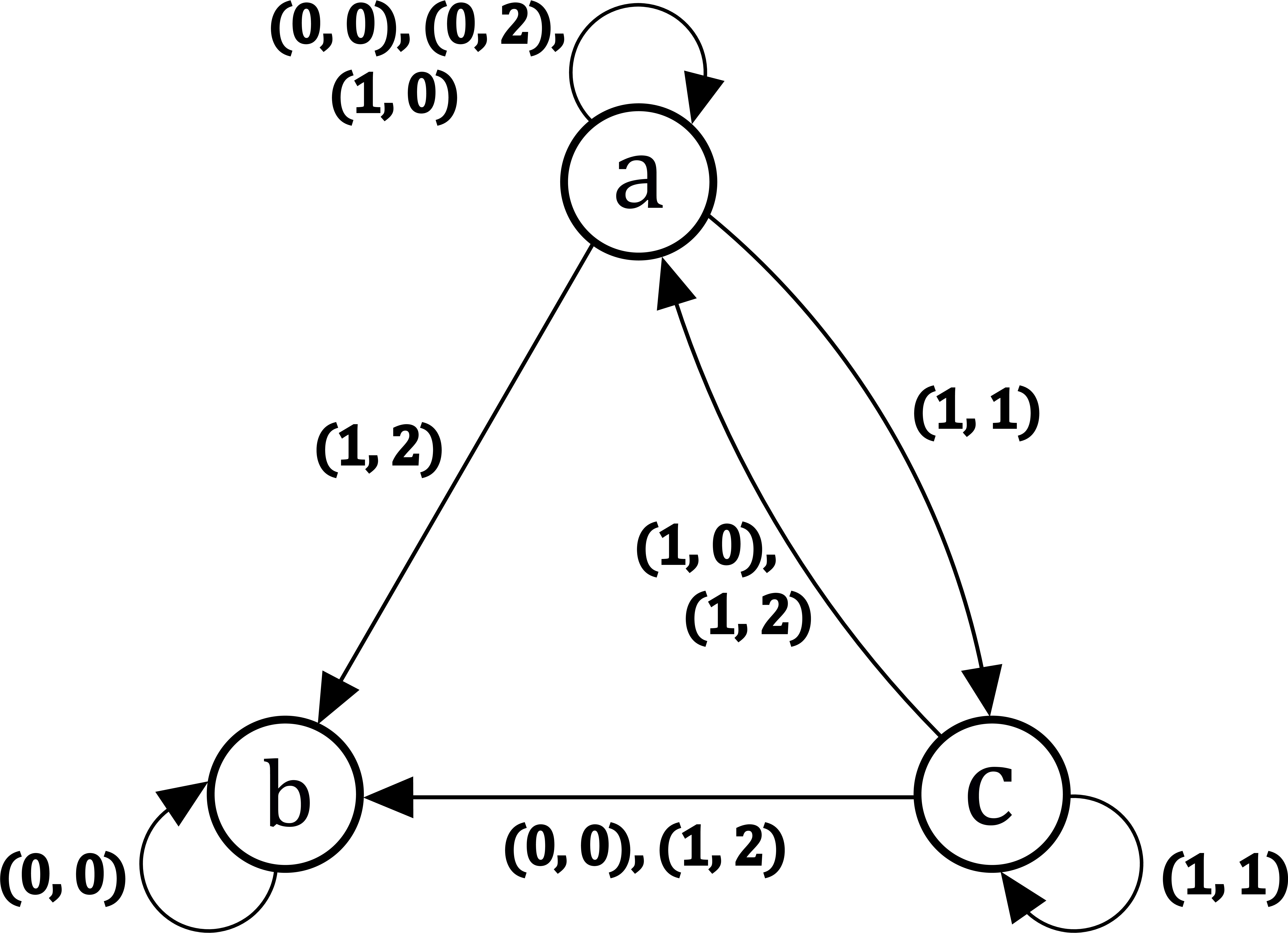}
\caption{A digraph labeled with $I = \{0, 1\} \times \{0, 1, 2\}$} \label{figure: digraph for planar 01}
\end{figure}

\noindent Then, the adjacency matrices are
\begin{equation*}
A_0 =
\begin{pmatrix}
2 & 0 & 0\\
0 & 1 & 0\\
0 & 1 & 0
\end{pmatrix},
\hspace{4pt} A_1 =
\begin{pmatrix}
1 & 1 & 1 \\
0 & 0 & 0 \\
2 & 1 & 1
\end{pmatrix}.
\end{equation*}
Let
\begin{align*}
C_{N, k} = 2^{N-k+1} \Big( ( 2 + 2\sqrt{2} )(1+ \sqrt{2})^k - ( 2 - 2\sqrt{2} ) (1- \sqrt{2})^k \Big).
\end{align*}
We have $\mathrm{dim}_{\mathrm{H}}(X) = \log_2{r} = 1.6416\cdots$, where $r$ satisfies
\begin{equation*}
r = \sum_{N =1}^\infty \left( \sum_{k = 0}^N
{C_{N, k}}^{\log_3{2}} \right) r^{-N-1} = 3.1201\cdots.
\end{equation*}
If we denote by $C_N$ the coefficient of $r^{-N}$, then $r$ is also the spectral radius of the operator
\begin{equation*}
\begin{pmatrix}
C_0 & C_1 & C_2 & \cdots \\
1 & 0 & 0 & \cdots \\
0 & 1 & 0 & \\
0 & 0 & 1  & \\
& \vdots & & \ddots \\
\end{pmatrix}.
\end{equation*}

\end{example*}

\spa

\begin{example*}

Let $G$ be a directed graph with $2$ vertices, and $I = \{0, 1, 2\} \times \{0, 1, 2, 3\}$. Consider a sofic system with the following adjacency matrices. (The digraph $G$ in this example has $19$ edges, so we omit it.)
\begin{equation*}
A_0 =
\begin{pmatrix}
1 & 0 \\
2 & 0 \\
\end{pmatrix},
\hspace{4pt} A_1 =
\begin{pmatrix}
2 & 1 \\
1 & 2 \\
\end{pmatrix},
\hspace{4pt} A_2 =
\begin{pmatrix}
3 & 2 \\
2 & 3 \\
\end{pmatrix}.
\end{equation*}
Then, $\mathrm{dim}_{\mathrm{H}}(X) = \log_3{r} = 1.6994\cdots$, where $r$ satisfies
\begin{align*}
r = \sum_{N=0}^{\infty} \left(\sum_{k = 0}^N
\begin{pmatrix}
N \\
k
\end{pmatrix}
\left( \frac{3^k \spa 5^{N-k} -1}{2} \right)^{\log_4{3}} \right) r^{-N} = 6.4693\cdots.
\end{align*}

\end{example*}

\subsection{Results in $\mathbb{R}^3$} \label{subsection: results in 3-dim} \hfill\\
\quad Next, we state our results on sofic sets in $\mathbb{R}^3$. Let $1 < m_1 \leq m_2 \leq m_3$ be natural numbers, and let $I = \{0, 1, \ldots, m_1-1\} \times \{0, 1, \ldots, m_2-1\} \times \{0, 1, \ldots, m_3-1\}$ be the set of indices used to label edges in a directed graph $G$. In a manner analogous to the planar case, we define a sofic set $X \subset \mathbb{R}^3$. Let $I_1 = \{0, 1, \ldots, m_1-1\}$ and $I_2 = \{0, 1, \ldots, m_2-1\}$. It has been proved that there is a matrix $A_{(s,t)} \in M_n(\mathbb{Z}_{\geq 0})$ for each $(s,t) \in I_1 \times I_2$ such that the Hausdorff dimension of $X$ is given by the following expression. (\cite[Theorems 1.1 and 4.1]{Barral--Feng}, \cite[Lemma 2.5]{Z.Feng})
\begin{align} \label{equation: dimension formula for 3d sofic}
\mathrm{dim_H}(X)
=
\lim_{N \to \infty} \spa \frac{1}{N} \log_{m_1}{\hspace{-3pt}\sum_{(s_1^{}, \ldots, s_N^{}) \in I_1^N} 
\left( \spa \sum_{ (t_1^{}, \ldots, t_N^{}) \in I_{2}^N} \hspace{-2pt}
{\norm{A_{(s_1^{},t_1^{})} \cdots A_{(s_N^{},t_N^{})}}
}^{a_{2}^{}}
\right)^{\hspace{-2pt} a_1^{}}}.
\end{align}
Here, $a_1 = \log_{m_2}{m_1}$ and $a_2 = \log_{m_3}{m_2}$.

We say that a sofic set $X \subset \mathbb{R}^3$ has a {\textbf{recursive structure}} with $\boldsymbol{v} \in \mathbb{R}^n$, if for any $s \in I_1$, there is $t \in I_2$ such that the image of $A_{(s, t)}$ is spanned by $\boldsymbol{v}^{\top}$. (Here, the matrices act on $\mathbb{R}^n$ from right.) Suppose $X$ satisfies this. For each $s \in I_1$, we will introduce (in \S \ref{subsection: tower decomposition for sofic sets in 3d}) a linear operator $M_s$ on a space indexed by a certain tree. The following theorem shows that we can reduce the double summations in equation (\ref{equation: dimension formula for 3d sofic}) to a single one (see Theorem \ref{theorem: tower decomposition for 3d} for the details).

\begin{theorem*}
Suppose that $X$ has a recursive structure, and that $\sum_{w \in I_1 \times I_2} A_w$ is primitive. Then,
\begin{equation*}
\mathrm{dim}_{\mathrm{H}}(X)
=
\lim_{N \to \infty} \spa \frac{1}{N} \log_{m_1^{}}{\hspace{-3pt}\sum_{(s_1^{}, \ldots, s_N^{}) \in I_1^N}
{\norm{ M_{s_N^{}} \spa M_{s_{N-1}^{}} \cdots \spa M_{s_1^{}} \spa \Phi_0 }_1
}^{a_1^{}}},
\end{equation*}
where $\Phi_0$ is a certain vector in the space indexed by a tree and $\norm{\cdot}_1$ is the $l^1$-norm.
\end{theorem*}

This expression can sometimes be further simplified (Proposition \ref{proposition: reapplication}), leading to explicit calculations of the dimension. The following examples illustrate such phenomena.

\begin{example*}

Let $I = \{0, 1\} \times \{0, 1, 2\} \times \{0, 1, 2, 3\}$. Consider the directed graph in Figure \ref{figure: digraph for 3d example} labeled with $I$.
\begin{figure}[h!]
\includegraphics[width=\textwidth-4cm]{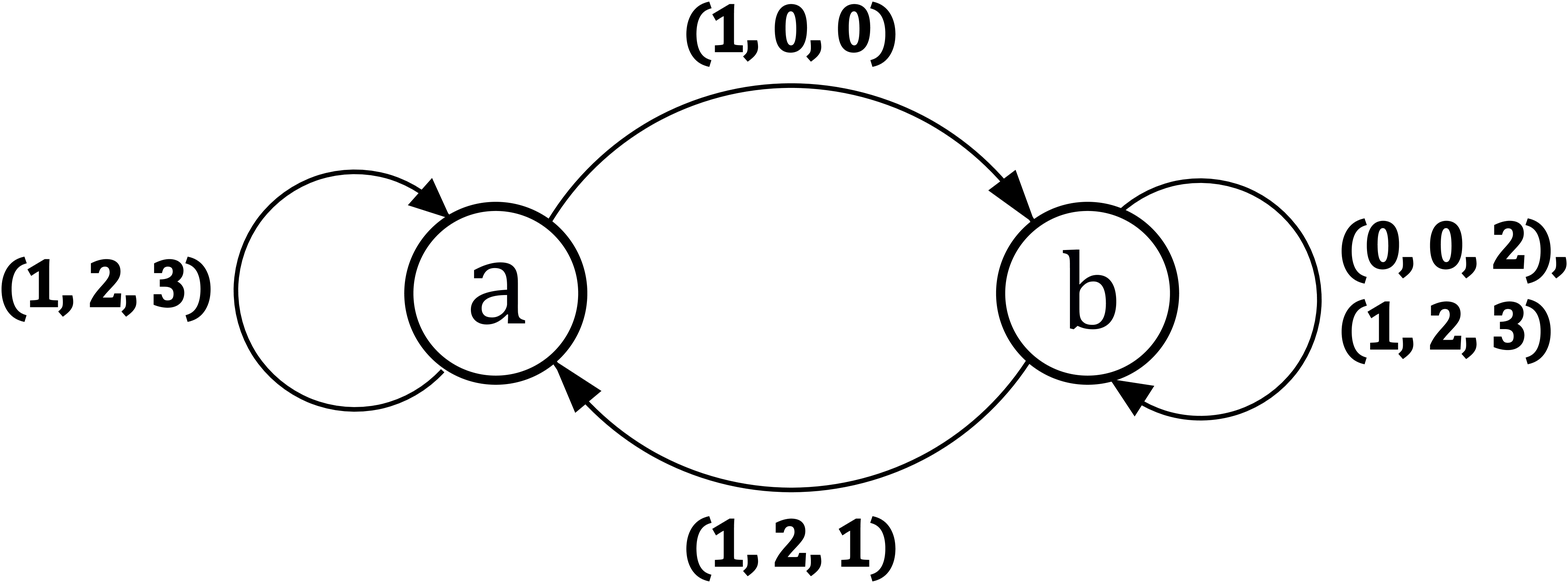}
\caption{A digraph labeled with $I = \{0, 1\} \times \{0, 1, 2\} \times \{0, 1, 2, 3\}$} \label{figure: digraph for 3d example}
\end{figure}
Then, the adjacency matrices are
\begin{equation*}
A_{(0, 0)} =
\begin{pmatrix}
0 & 0 \\
0 & 1 \\
\end{pmatrix},
\hspace{4pt} A_{(1, 0)} =
\begin{pmatrix}
0 & 1 \\
0 & 0 \\
\end{pmatrix},
\hspace{4pt} A_{(1, 2)} =
\begin{pmatrix}
1 & 0 \\
1 & 1 \\
\end{pmatrix}.
\end{equation*}
For a natural number $k$, let
\begin{equation*}
b_k =
\norm{
{\begin{pmatrix}
0 & 1 & 2^{\log_4{3}} & 3^{\log_4{3}} & 4^{\log_4{3}} & \cdots \\
1 & 0 & 0 & \cdots \\
0 & 1 & 0 & 0 & \cdots \\
0 & 0 & 1 & 0 & 0 & \cdots \\
& \vdots & & & \ddots \\
\end{pmatrix}}^k
\begin{pmatrix}
1 \\
0 \\
0 \\
0 \\
\vdots \\
\end{pmatrix}
}_{1}^{\log_3{2}}.
\end{equation*}
Then, $\mathrm{dim}_{\mathrm{H}}(X) = \log_2{r} = 1.1950\cdots$, where $r$ satisfies
\begin{align*}
r
= \sum_{k = 0}^{\infty} \frac{b_k}{r^k}
&= 1 + \frac{1}{r} + \frac{\sqrt{2}}{r^2} + \frac{ (2 + 2^{\log_4{3}})^{\log_3{2}} }{r^3} + \frac{ (3 + 2^{\log_4{3}} + 3^{\log_4{3}})^{\log_3{2}} }{r^4} + \cdots & \\
&  = 2.2894\cdots. &
\end{align*}

\end{example*}

\begin{example*}

Let $I = \{0, 1, 2\} \times \{0, 1, 2, 3\} \times \{0, 1, 2, 3, 4\}$. Consider the directed graph in Figure \ref{figure: digraph for 3d example 2} labeled with $I$.
\begin{figure}[h!]
\includegraphics[width=\textwidth-4cm]{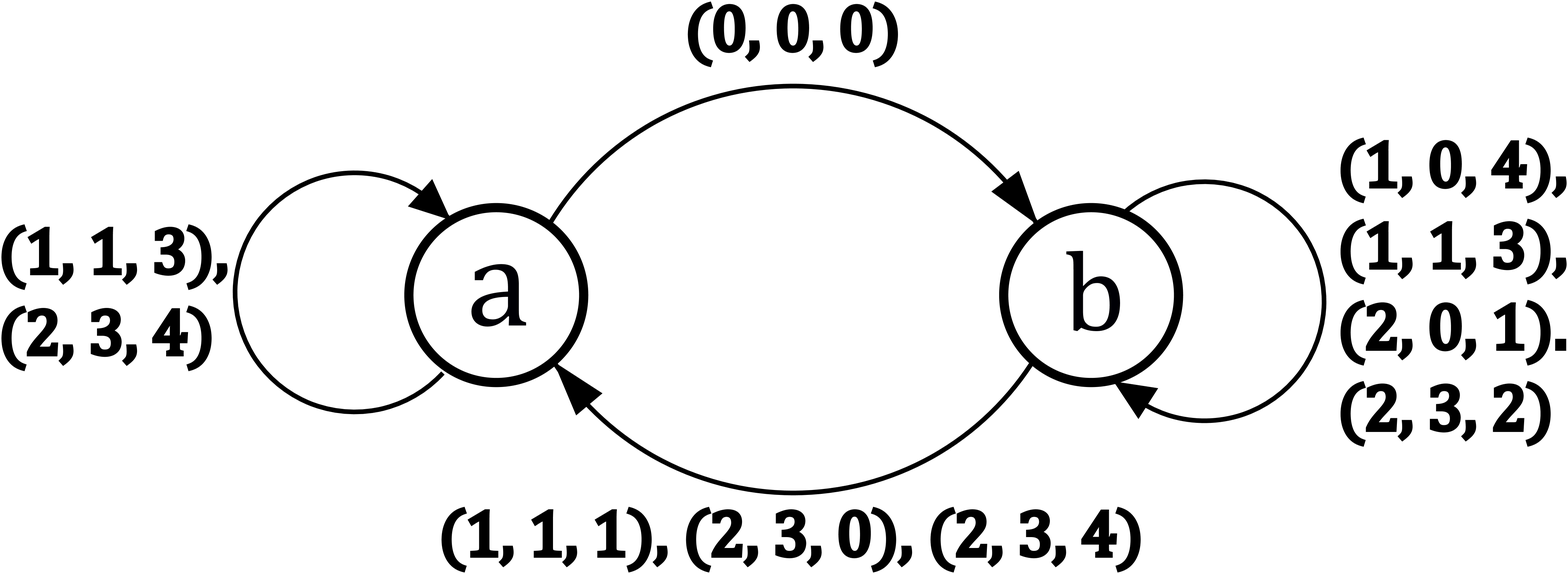}
\vspace{3pt}%
\caption{A digraph labeled with $I = \{0, 1, 2\} \times \{0, 1, 2, 3\} \times \{0, 1, 2, 3, 4\}$} \label{figure: digraph for 3d example 2}
\end{figure} \\
Then, the adjacency matrices are
\begin{equation*}
A_{(0, 0)} =
\begin{pmatrix}
0 & 1 \\
0 & 0 \\
\end{pmatrix},
\hspace{4pt} A_{(1, 0)} =
\begin{pmatrix}
0 & 0 \\
0 & 1 \\
\end{pmatrix},
\hspace{4pt} A_{(1, 1)} =
\begin{pmatrix}
1 & 0 \\
1 & 1 \\
\end{pmatrix},
\hspace{4pt} A_{(2, 0)} =
\begin{pmatrix}
0 & 0 \\
0 & 1 \\
\end{pmatrix},
\hspace{4pt} A_{(2, 3)} =
\begin{pmatrix}
1 & 0 \\
2 & 1 \\
\end{pmatrix}.
\end{equation*}
Let $a_1^{} = \log_4{3}$ and $a_2^{} = \log_5{4}$. Define
\begin{align*}
\hspace{10pt} b_N
&= \sum_{(s_1, \ldots, s_N) \in \{1, 2\}^N} \Bigg( \big( N + \#\left\{ j \setcond 1 \leq j \leq N, s_j = 2\right\} \big)^{a_2^{}} \hspace{20pt} &\\
&\hspace{150pt} + \sum_{k = 1}^{N-1} 2^{k-1} \big( N + \#\left\{j \setcond k \leq j \leq N \text{\hspace{1pt} and \hspace{1pt}} s_j = 2\right\} \big)^{a_2^{}} \Bigg)^{a_1^{}}. &
\end{align*}
Then, $\mathrm{dim}_{\mathrm{H}}(X) = \log_2{r} = 2.224\cdots$, where $r$ satisfies
\begin{align*}
r
= \sum_{k = 0}^{\infty} \frac{b_k}{r^k}
= 4.673\cdots.
\end{align*}

\end{example*}

\begin{remark*}
Regarding the dimension theory of sofic sets, \cite{Olivier} discussed the uniqueness of the measure with full dimension in the planar case. Additionally, \cite{Z.Feng} studied the conditions under which general self-affine fractals, including sofic sets, have the same Hausdorff dimension and box dimension.
\end{remark*}

\section{Background} \label{section: background}

This section introduces the definition of sofic systems and sofic sets, and then proceeds to their dimension formula in terms of matrix products. Let $\mathbb{T} = \mathbb{R}/\mathbb{Z}$. We begin by reviewing the definition of sofic systems. Weiss \cite{Weiss} defined {\it{sofic systems}} as subshifts which are factors of shifts of finite type. Using the results by \cite{Boyle--Kitchens--Marcus}, sofic systems can also be characterized in the following way.

\begin{definition}[{{\cite[Proposition 3.6]{Kenyon--Peres: sofic}}}] \label{definition: sofic systems}
Consider a finite directed graph $G = \langle V, E \rangle$ with possible loops and multiple edges. Let $I$ be a set of labels. Suppose that edges of $G$ are labeled using elements of $I$ in a ``right-resolving'' fashion: every two edges emanating from the same vertex have different labels. Then, the set $S \subset I^{\mathbb{N}}$ of sequences of labels that arise from infinite paths in $G$ is called the \textbf{sofic system}.
\end{definition}

Fix a positive integer $d$, the Euclidean dimension of the space in which we will be working. Consider natural numbers $1 < m_1 \leq m_2 \leq \cdots \leq m_d$, and define $I_i = \{0, 1, \ldots, m_i-1\}$ for $1 \leq i \leq d$. Let $I = I_1 \times I_2 \times \cdots \times I_d$, and define $R: I^{\mathbb{N}} \rightarrow {\mathbb{T}}^d$ by
\begin{equation*}
R \left( \left( e^{(n)}_1, e^{(n)}_2, \ldots, e^{(n)}_d \right)_{n=1}^{\infty} \right) = \left( \sum_{n=1}^{\infty} \frac{e^{(n)}_1}{{m_1}^n}, \cdots, \sum_{n=1}^{\infty} \frac{e^{(n)}_d}{{m_d}^n} \right).
\end{equation*}

Let $S$ be a sofic system with the set of labels $I$. Then, $R(S)$ is a compact set invariant under the map $T_A: {\mathbb{T}}^d \rightarrow {\mathbb{T}}^d$, which is induced from the multiplication of the diagonal matrix $A = \mathrm{diag}(m_1, m_2, \ldots, m_d)$. This set $R(S)$ is referred to as the \textbf{sofic set}. When $d=2$, Kenyon and Peres \cite{Kenyon--Peres: sofic} proved that the Hausdorff dimension of sofic sets can be expressed as the limit of a certain combinatorial sum involving matrix products. This result extends to general $d$, as we shall explain.

Suppose $S$ is a sofic system defined from a directed graph $G = \langle V, E \rangle$ with $n$ number of vertices; $V = \{v_1, v_2, \ldots, v_n\}$. Let $W_i = I_1 \times \cdots \times I_i$ for $1 \leq i \leq d-1$, and for each $s \in W_{d-1}$, define the restricted adjacency matrix $A_{s} = \left( a_{ij}(s) \right)_{i,j} \in \mathrm{M}_n(\mathbb{Z}_{\geq 0})$ by
\[ a_{ij}(s) = \# \left\{ 
				e \in E \setcond \text{$e$ is from vertex $v_i$ to $v_j$ and its label is $(s, k) \in I$ for some $k \in I_d$} 
			\right\}. \]
For a natural number $N$ and $s = (s_1^{}, \ldots, s_N^{}) \in W_i^N$, define
\[ W_{i+1}^N(s) = \left\{ s' \in W_{i+1}^N \setcond \text{ There is $(t_1^{}, \ldots, t_N^{}) \in I_{i+1}^N$ with $s' = \big( (s_1^{}, t_1^{}), \ldots, (s_N^{}, t_N^{}) \big)$  } \right\}. \]

The following theorem is a special case of \cite[Theorems 1.1 and 4.1]{Barral--Feng}, as explained in \cite[Lemma 2.5]{Z.Feng}.

\begin{theorem}[{{\cite[Lemma 2.5]{Z.Feng}}}] \label{theorem: combinatorial formula for the Hausdorff dimension of sofic sets}
Let $a_i^{} = \log_{m_{i+1}^{}}{m_i^{}} \hspace{2pt}$ for $ \spa 1 \leq i \leq d-1$. For any sofic system $S \subset I^{\mathbb{N}}$, we have
\begin{equation*} 
\mathrm{dim}_{\mathrm{H}}(R(S))
= \lim_{N \to \infty} \spa \frac{1}{N} \log_{m_1^{}}{\hspace{-5pt}\sum_{s_1^{} \in W_1^N} \hspace{-3pt}
\left( \hspace{-1pt} \sum_{s_2^{} \in W_2^N(s_1^{})} \hspace{-3pt}
	\left( \hspace{-2pt} \cdots \hspace{-2pt}
		\left( \sum_{\substack{
					s_{d-1}^{} \in W_{d-1}^N(s_{d-2}^{}) \\
					s_{d-1}^{} = (s^{(1)}, \ldots, s^{(N)})}} \hspace{-12pt}
		{\norm{ A_{s^{(1)}} \cdots A_{s^{(N)}} } }^{a_{d-1}^{}} \hspace{-2pt}
		\right)^{\hspace{-2pt} a_{d-2}^{}} \hspace{-14pt} \cdots
	\right)^{\hspace{-2pt} a_2^{}} \hspace{1pt}
\right)^{\hspace{-2pt} a_1^{}}}.
\end{equation*}
Here, the norm for matrices can be arbitrary.
\end{theorem}

This theorem can also be proved using a concept called {\it{weighted topological entropy}} as discussed in \cite[Claim 1.6]{Alibabaei}. In \cite{Alibabaei}, the Hausdorff dimension of a sofic set in $\mathbb{T}^3$ is computed, although the example considered there was trivial in the following sense. Suppose that either all the matrices commute and $\sum_{u \in W_{r-1}} A_u$ is primitive, or the matrices share a common eigenvector. In either case, it can be shown that
\begin{equation*}
\mathrm{dim}_{\mathrm{H}}(R(S)) = \log_{m_1^{}}{\sum_{s_1^{} \in W_1} 
\left(  \sum_{s_2^{} \in W_2(s_1^{})} 
	\left(  \cdots
		\left( \sum_{s_{d-1}^{} \in W_{d-1}(s_{d-2}^{})} 
			{\lambda_{\spa s_{d-1}^{}}}^{a_{d-1}^{}}
		\right)^{a_{d-2}^{}} \cdots
	\right)^{ a_2^{}} 
\right)^{a_1^{}}},
\end{equation*}
where $\lambda_{\spa u}$ is the spectral radius of $A_u$. The aim of this paper is to provide sophisticated calculations of the dimension for much more non-trivial cases.

For planar cases, the following structure behind the dimension formula is explained after Proposition 3.4 in \cite{Kenyon--Peres: sofic}. Let $\flo{L}$ be an operator acting on the continuous functions on $S^{\spa \# V -1}$ by
\[ (\flo{L}f)(x) = \sum_{s = 0}^{m_1^{} -1} \norm{A_s x}^{a_1^{}}
f \left( \frac{A_s \spa x}{ \norm{A_s \spa x} } \right). \]
Then, we have $\mathrm{dim_H} (X) = \log_{m_1^{}} {\rho(\flo{L})}$, where $\rho(\flo{L})$ is the spectral radius of $\flo{L}$. When one of the matrices has a $1$-dimensional image, the sphere $S^{\spa \#V -1}$ can be replaced with a countable set of directions, resulting in an expression for $\flo{L}$ as an infinite matrix. This structure provides valuable intuition for the arguments in \S \ref{subsection: tower decomposition for planar sofic sets}. In contrast, we were not able to uncover this type of ``bigger picture'' in $3$-dimensional cases due to the nesting of the summations. At present, the technique presented in \S \ref{subsection: tower decomposition for sofic sets in 3d} appears to be a ``fortunate discovery''. Hopefully, a deeper underlying structure will be uncovered in the future.

\section{Theorems and proofs} \label{section: theorems and proofs}

\subsection{Tower decomposition for planar sofic sets} \label{subsection: tower decomposition for planar sofic sets} \hfill\\
\quad In this section, we generalize the technique introduced by Kenyon and Peres \cite[Example 4.2]{Kenyon--Peres: sofic} to compute the Hausdorff dimension of sofic sets in $\mathbb{T}^2$. Recall that $m_1 \leq m_2$ are natural numbers and $I_1 = \{0, 1, \ldots, m_1-1\}$. Let $X \subset \mathbb{T}^2$ be a sofic set introduced in \S \ref{section: background}. By Theorem \ref{theorem: combinatorial formula for the Hausdorff dimension of sofic sets}, the Hausdorff dimension of $X$ is given by the following expression, where $\alpha = \log_{m_2}{m_1}$.
\begin{equation*}
\mathrm{dim}_{\mathrm{H}}(X)
=
\lim_{N \to \infty} \spa \frac{1}{N} \log_{m_1^{}}{ \sum_{
				(u_1^{}, \ldots, u_N^{}) \in I_1^{N}}
			{\norm{A_{u_1^{}} \cdots A_{u_N^{}}}
			}^{\alpha}}.
\end{equation*}
We define the norm as the sum of the absolute values of all entries of the matrix.

\begin{breakdefinition}
(1) A matrix $A$ is said to be \textbf{primitive} if there is an integer $d$ such that every entry in $A^d$ is positive.

\noindent(2) For $A \in M_n(\mathbb{R})$, define $\mathrm{Im_R^{}}(A) = \left\{ x^{\top} A \setcond x \in \mathbb{R}^n \right\}$.
\end{breakdefinition}

For $s \in I_1^L$ with a natural number $L$, define the \textbf{exclusion set} for each natural number $k$ by
\[ I_1^k \backslash \{s\} = \left\{ u \in I_1^k \setcond \text{The string $s$ does not appear in $u$} \right\}. \]
Finally, let $A = \sum_{i \in I_1} A_i$. We are now ready to state and prove the following main theorem for planar sofic sets.

\begin{theorem} \label{theorem: tower decomposition for 2d}
Suppose $A$ is primitive. Additionally, assume we have an integer $L > 0$ and a string $s = (s_1^{}, \ldots, s_L^{}) \in I_1^L$ such that there is $\boldsymbol{v} \in \mathbb{R}^n \backslash \{0\}$ with
\[ \mathrm{Im_R^{}}( A_{s_1^{}} \cdots A_{s_L^{}} ) = \mathrm{Span}\{ \boldsymbol{v}^\top \}. \]
For $u = (u_1^{}, \ldots, u_k^{}) \in I_1^k$, define $J_u \geq 0$ by
\[ \boldsymbol{v}^{\top} \spa A_{u_1^{}} \cdots A_{u_{k}^{}} \hspace{2pt} A_{s_1^{}} \cdots A_{s_L^{}} = J_u \spa {\boldsymbol{v}}^{\top}. \]
Set $C_k = \sum_{u \in I_1^k \backslash \{s\}} J_u^{\alpha}$ for each non-negative integer $k$. Then, we have
\[\mathrm{dim}_{\mathrm{H}}(X) = \log_{m_1^{}}{r}, \]
where $r$ is the unique positive solution to the equation
\[ r^L = C_0 + \frac{C_1}{r} + \frac{C_2}{r^2} + \cdots. \]
\end{theorem}

\begin{proof}
We break down the summand by constructing what can be described as a tower-decomposition. Specifically, we define a vector $\Phi_N = \left( \Phi_{N, k} \right)_{k = 0}^{\infty}$ for each natural number $N \geq L$, whose norm closely approximates that of the original matrix product. Let $\mathbb{e}$ be the column vector with $1$ in every entry. Let
\[ \Phi_{N, N} = \left[ \boldsymbol{v}^{\top} \mathbb{e} \right]^{\alpha}, \]
\[ \Phi_{N, N-L} = \left[ \boldsymbol{v}^{\top} \spa A_{s_1^{}} \cdots A_{s_{L}^{}} \spa \mathbb{e} \right]^{\alpha}, \]
and when $0 \leq k < N-L$, let
\[
\Phi_{N, k}
= \sum_{(u_i)_i \in I_1^{N-L-k}} \left[ \boldsymbol{v}^{\top} A_{u_{1}^{}} \ldots A_{u_{N-L-k}^{}} \spa A_{s_1^{}} \cdots A_{s_{L}^{}} \spa \mathbb{e} \right]^{\alpha}.
\]
Set $\Phi_{N, k} = 0$ for all other $k$.

\begin{claim} \label{claim: tower decomposition}
We have the following relation:
\begin{equation*}
\mathrm{dim}_{\mathrm{H}}(X)
=
\lim_{N \to \infty} \frac{1}{N} \log_{m_1^{}} \norm{ \Phi_N }_1,
\end{equation*}
where $\norm{\cdot}_1$ is the $l^1$-norm.
\end{claim}
\begin{proof}[Proof of Claim \ref{claim: tower decomposition}]
Since $A = \sum_{i \in I_1} A_i \in \mathrm{M}_n(\mathbb{Z}_{\geq0})$ is primitive, there is an integer $d$ such that every entry of $A^d$ is at least 1. From this, we can conclude that
\[ \norm{ \boldsymbol{x}^{\top} }_1
\leq
\norm{ \boldsymbol{x}^{\top} \spa A^{d} \spa A_{s_1^{}} \cdots A_{s_L^{}} }_1 \]
for every $\boldsymbol{x} \in (\mathbb{R}_{\geq 0})^n$. Also, for each natural number $k$, we can take a string $t^{(k)} = (t_1^{(k)}, \ldots, t_k^{(k)}) \in I_1^k$ such that $A_{t_1^{(k)}} \cdots A_{t_k^{(k)}} \ne O$, again because $A$ is primitive. Then, we have
\[ \norm{ \boldsymbol{x}^{\top} }_1
\leq
\norm{ \boldsymbol{x}^{\top} \spa A^{d} \spa A_{t_1^{(k)}} \cdots A_{t_k^{(k)}} }_1. \]
Using the fact that $(x + y)^{\alpha} \leq x^{\alpha} + y^{\alpha}$ for $x, y \geq 0$ when $0 \leq \alpha \leq 1$, we can make the following evaluation:
\begin{align*}
\norm{ \Phi_N }_1
& \leq
\sum_{k = 0}^{N-L} \hspace{4pt} \sum_{(u_i)_i \in I_1^{N-L-k}}
\left[ \boldsymbol{v}^{\top} A_{u_1^{}} \cdots A_{u_{N-L-k}^{}} \spa A_{s_1^{}} \cdots A_{s_L^{}} \spa A^d \spa A_{t_1^{(k)}} \cdots A_{t_k^{(k)}} \spa \mathbb{e} \right]^{\alpha} + \Phi_{N, N} \\
& =
\sum_{k = 0}^{N-L} \hspace{4pt} \sum_{(u_i)_i \in I_1^{N-L-k}}
\left[ \boldsymbol{v}^{\top} A_{u_1^{}} \cdots A_{u_{N-L-k}^{}} \spa A_{s_1^{}} \cdots A_{s_L^{}} \spa \left( \sum_{q = 0}^{m-1} A_q \right)^d \spa A_{t_1^{(k)}} \cdots A_{t_k^{(k)}} \spa \mathbb{e} \right]^{\alpha}
+ \left[ \boldsymbol{v}^{\top} \mathbb{e} \right]^{\alpha} \\
& \leq
m_1^d \spa N \spa \sum_{(u_i^{})_i \in I_1^{N + d}}
		{\left[ \boldsymbol{v}^{\top} A_{u_{1}^{}} \cdots A_{u_{N+d}^{}} \mathbb{e}
		\right]
		}^{\alpha}
+ \left[ \boldsymbol{v}^{\top} \mathbb{e} \right]^{\alpha} \\
& \leq
m_1^d \spa N \spa \norm{\boldsymbol{v}}_{\infty} \sum_{(u_i^{})_i \in I_1^{N + d}}
		{\left[ \mathbb{e}^{\top} A_{u_{1}^{}} \cdots A_{u_{N+d}^{}} \mathbb{e}
		\right]
		}^{\alpha}
+ \left[ \boldsymbol{v}^{\top} \mathbb{e} \right]^{\alpha} \\
& \leq
m_1^d \spa N \spa \norm{\boldsymbol{v}}_{\infty} \sum_{(u_i^{})_i \in I_1^{N + d}}
		{\left[ \boldsymbol{v}^{\top} \spa A^{d} \spa A_{u_{1}^{}} \cdots A_{u_{N+d}^{}} \spa A^d \spa A_{s_1^{}} \cdots A_{s_L^{}} \spa \mathbb{e}
		\right]
		}^{\alpha}
+ \left[ \boldsymbol{v}^{\top} \mathbb{e} \right]^{\alpha} \\
& \leq
m_1^{3d} \spa N \spa \norm{\boldsymbol{v}}_{\infty} \spa \norm{\Phi_{N+3d+L}}_1.
\end{align*}
Taking the logarithm and letting $N \rightarrow \infty$, we have
\begin{equation*}
\lim_{N \to \infty} \frac{1}{N} \log_{m_1^{}} \norm{ \Phi_N }_1
= \lim_{N \to \infty} \spa \frac{1}{N} \log_{m_1^{}}{ \sum_{
				(u_1^{}, \ldots, u_N^{}) \in I_1^{N}}
			{\norm{A_{u_1^{}} \cdots A_{u_N^{}}}
			}^{\alpha}}.
\end{equation*}
\end{proof}

Now, the vector $\boldsymbol{v}^{\top} A_{u_{1}^{}} \ldots A_{u_{N-L-k}^{}} \spa A_{s_1^{}} \cdots A_{s_{L}^{}}$ in the definition of $\Phi_{N,k}$ is already a constant multiple of $\boldsymbol{v}^{\top}$. Therefore, for any $L-1 \leq k \leq N-L$ we have
\begin{flalign*}
&\; & \sum_{(w_1^{}, \ldots, w_{k-L+1}^{}) \in I_1^{k-L+1} \backslash \{s\}} \hspace{2pt} \sum_{(u_i)_i \in I_1^{N-L-k}} \left[ \boldsymbol{v}^{\top} \spa A_{u_{1}^{}} \ldots A_{u_{N-L-k}^{}} \spa A_{s_1^{}} \cdots A_{s_{L}^{}} \spa A_{w_1^{}} \cdots A_{w_{k-L+1}^{}} \spa A_{s_1^{}} \cdots A_{s_{L}^{}} \spa \mathbb{e} \right]^{\alpha} &
\end{flalign*}
\begin{flalign*}
& \hspace{15pt} = \sum_{w \in I_1^{k-L+1} \backslash \{s\}} J_w^{\alpha} \Phi_{N, k} & \\
& \hspace{15pt} = C_{k-L+1} \Phi_{N, k}.&
\end{flalign*}
Also,
\begin{flalign*}
&\sum_{(w_1^{}, \ldots, w_{N-L+1}^{}) \in I_1^{N-L+1} \backslash \{s\}} \hspace{2pt} \left[ \boldsymbol{v}^{\top} \spa A_{w_1^{}} \cdots A_{w_{N-L+1}^{}} \spa A_{s_1^{}} \cdots A_{s_{L}^{}} \spa \mathbb{e} \right]^{\alpha}
= \sum_{w \in I_1^{N-L+1} \backslash \{s\}} J_w^{\alpha} \Phi_{N, N} &
\end{flalign*}
\begin{flalign*}
& \hspace{273.5pt}
= C_{N-L+1} \Phi_{N, N}. &
\end{flalign*}

We thus conclude that
\[
\Phi_{N+1, 0} = \sum_{k = 0}^{\infty} \spa C_{k} \spa \Phi_{N, k+L-1}.
\]
Also, we have $\Phi_{N+1, k+1} = \Phi_{N, k}$. Thus, we can express the sum in question using a linear operator as follows. Let $\mathbb{R}^{\bigoplus \mathbb{N}_0}$ denote the direct sum of $\mathbb{R}$ indexed by $\mathbb{N}_0$; an element in this direct sum has only finitely many non-zero entries. Define $M: \mathbb{R}^{\bigoplus \mathbb{N}_0} \rightarrow \mathbb{R}^{\bigoplus \mathbb{N}_0}$ by \\[-6pt]
\begin{equation*}
\begin{matrix}
\phantom{
\begin{matrix}
\overbrace{
    \hphantom{\begin{matrix}
M \end{matrix}}
}^{\text{$L$}}
\end{matrix}
\phantom{
\begin{matrix} M \end{matrix}}}\\
M =
\end{matrix}
\begin{matrix}
\begin{matrix}
\overbrace{
    \hphantom{\begin{matrix}
0 & 0 & 0 & \cdots & 0 \end{matrix}}
}^{\text{$(L-1)$ columns}}
\end{matrix}
\phantom{
\begin{matrix} 0 & C_0 & C_1 & C_2 & \cdots \end{matrix}}\\
\begin{pmatrix}
0 & 0 & 0 & \cdots & 0 & C_0 & C_1 & C_2 & \cdots \\
1 & 0 & 0 & \cdots &    &  0   &  0   &  \cdots \\
0 & 1 & 0 & 0 & \cdots \\
0 & 0 & 1 & 0 & \cdots \\
& \vdots & & \ddots \\
\end{pmatrix}.
\end{matrix}
\end{equation*} \\[4pt]
Then, $\Phi_{N+1} = M \Phi_N$. This yields for $N > L$,
\[
\norm{\Phi_N}_1 = \norm{M^{N-L} \spa \Phi_L }_1.
\]
Thus,
\[
\mathrm{dim}_{\mathrm{H}}(X)
=
\lim_{N \to \infty} \log_{m_1^{}}{ \norm{M^{N-L} \spa \Phi_L}_1^{\frac{1}{N}} }.
\]

Let $r$ be the unique solution to the following equation (which exists since $C_k$ grows at most exponentially).
\begin{equation*}
r = \sum_{k = 0}^{\infty} C_k r^{-k-L+1}.
\end{equation*}
Define $\boldsymbol{b} \in \mathbb{R}^{\mathbb{N}}$ by $\boldsymbol{b} = (1, r^{-1}, r^{-2}, \spa \ldots)^{\top}$. We see that $M \boldsymbol{b}$ is well-defined and $M \boldsymbol{b} = r\boldsymbol{b}$, which implies
\[
\norm{ M^{N-L} \spa \Phi_L }_1 \leq r^L \spa \norm{\Phi_L}_1 \spa \norm{ M^{N-L} \spa \boldsymbol{b} }_1 = \norm{\Phi_L}_1 \spa \norm{ \boldsymbol{b} }_1 \spa r^N.
\]
The other direction is done via an application of Perron-Frobenius theorem. Let $M_k$ be the upper $k \times k$ submatrix of $M$. By the Perron-Frobenius theorem,
\[
\liminf_{N \to \infty} \norm{ M^{N-L} \spa \Phi_L }_1^{\frac{1}{N}}
\geq \liminf_{N \to \infty} \norm{ M_k^{N-L} \spa \Phi_L^{(k)} }_1^{\frac{1}{N}} = r_k.
\]
Here, $r_k$ is the spectral radius of $M_k$ which satisfies
\begin{equation*}
r_k = \sum_{j = 0}^{k} C_j {r_k}^{-j-L+1},
\end{equation*}
Since $r_k \to r$ as $k \to \infty$, we obtain
\[
\mathrm{dim}_{\mathrm{H}}(X)
= \log_{m_1^{}}{r}. \qedhere
\]
\end{proof}

\begin{remark}
For comparison with the operator appearing in the next section \S\ref{subsection: tower decomposition for sofic sets in 3d} (Figure \ref{figure: operator in 3d}), we look at the operator $M$ as the sum of a shift operator and a ``return'' map. See Figure \ref{figure: operator in 2d}.
\begin{figure}[h!]
\includegraphics[width=\textwidth-2.5cm]{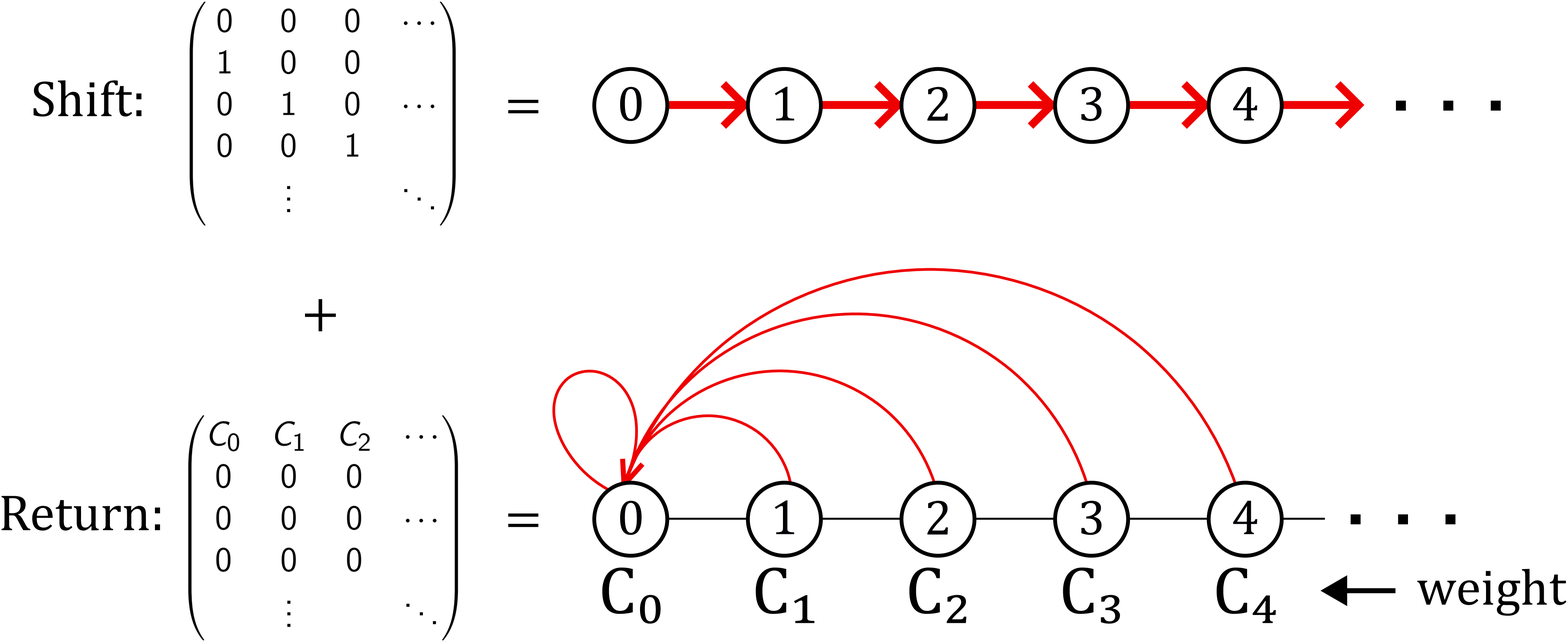}
\vspace{3pt}%
\caption{A description of the operator $M$ when $L = 1$} \label{figure: operator in 2d}
\end{figure}
\end{remark}

\subsection{Tower decomposition for sofic sets in $\mathbb{R}^3$} \label{subsection: tower decomposition for sofic sets in 3d} \hfill \\
\quad Let $m_1 \leq m_2 \leq m_3$ be natural numbers, and consider a sofic set $X \subset \mathbb{T}^3$ introduced in \S \ref{section: background}. Let $I_1 = \{0, 1, \ldots, m_1-1\}$ and $I_2 = \{0, 1, \ldots, m_2-1\}$.

By Theorem \ref{theorem: combinatorial formula for the Hausdorff dimension of sofic sets}, for each pair $(s,t) \in I_1 \times I_2$, there is a matrix $A_{(s,t)} \in M_n(\mathbb{Z}_{\geq 0})$ such that the Hausdorff dimension of $X$ can be expressed as
\begin{align*}
\mathrm{dim}_{\mathrm{H}}(X)
=
\lim_{N \to \infty} \spa \frac{1}{N} \log_{m_1^{}}{\hspace{-3pt}\sum_{(s_1^{}, \ldots, s_N^{}) \in I_1^N} 
\left( \spa \sum_{ (t_1^{}, \ldots, t_N^{}) \in I_{2}^N} \hspace{-2pt}
{\norm{A_{(s_1^{},t_1^{})} \cdots A_{(s_N^{},t_N^{})}}
}^{a_{2}^{}}
\right)^{\hspace{-2pt} a_1^{}}}.
\end{align*}
Here, $a_1^{} = \log_{m_2^{}}{m_1^{}}$ and $a_2^{} = \log_{m_3^{}}{m_2^{}}$.

\begin{definition}
We say that a sofic set $X \subset \mathbb{T}^3$ has a {\textbf{recursive structure}} with $\boldsymbol{v} \in \mathbb{R}^n$ if for any $s \in I_1$ there is $t \in I_2$ with $\mathrm{Im_R^{}}( A_{(s,t)} ) = \mathrm{Span} \{ \boldsymbol{v}^{\top} \}$. When this condition is satisfied, we define for each $u \in I_1$
\[
P(u) = \left\{ t \in I_2 \setcond \text{ Either $A_{(u, t)} = O$ or $\mathrm{Im_R^{}}( A_{(u, t)} ) = \mathrm{Span} \{ \boldsymbol{v}^{\top} \}$} \right\}.
\]
We say that $u \in I_1$ is \textbf{removable} if $P(u) = I_2$. 
\end{definition}

In the rest of this subsection, we assume that a sofic set $X$ has a recursive structure with $\boldsymbol{v}$. Define $J$ to be the collection of $u \in I_1$ that are not removable. As before, for $s \in I_1^N$ let $W_2^N(s) = \left\{ w \in (I_1 \times I_2)^N \setcond \text{The projection of $w$ onto $I_1^N$ is $s$.} \right\}$. Also, let $Q(u) = \{ u \} \times P(u)$ for $u \in I_1$.

Take $u \in I_1$. Define the constant $D_w(p) \geq 0$ for $w = (w_1^{}, \ldots, w_N^{}) \in (I_1 \times I_2)^N$ and $p \in Q(u)$ by
\[ \boldsymbol{v}^{\top} \spa A_{w_1^{}} \cdots A_{w_N^{}} \spa A_p = D_w(p) \boldsymbol{v}^{\top}, \]
and $D_{\varnothing}(p)$ by $\boldsymbol{v}^{\top} \spa A_p = D_{\varnothing}(p) \boldsymbol{v}^{\top}.$ Define $C_{\varnothing}(u) \geq 0$ by 
\[ C_{\varnothing}(u) = \sum_{p \in Q(u)} {D_{\varnothing}(p)}^{a_2}, \]
and $C_s(u) \geq 0$ for $s \in J^N$ by
\[ C_s(u) = \sum_{p \in Q(u)} \spa \spa \sum_{w \in R^N(s)} {D_w(p)}^{a_2}, \]
where $R^N(s) = \left\{ (s_i^{}, t_i^{})_i \in (I_1 \times I_2)^N \setcond \text{$t_i^{} \notin P(s_i^{})$ for each $1 \leq i \leq N$} \right\}$.

\begin{figure}[h!]
\includegraphics[width=\textwidth-1.5cm]{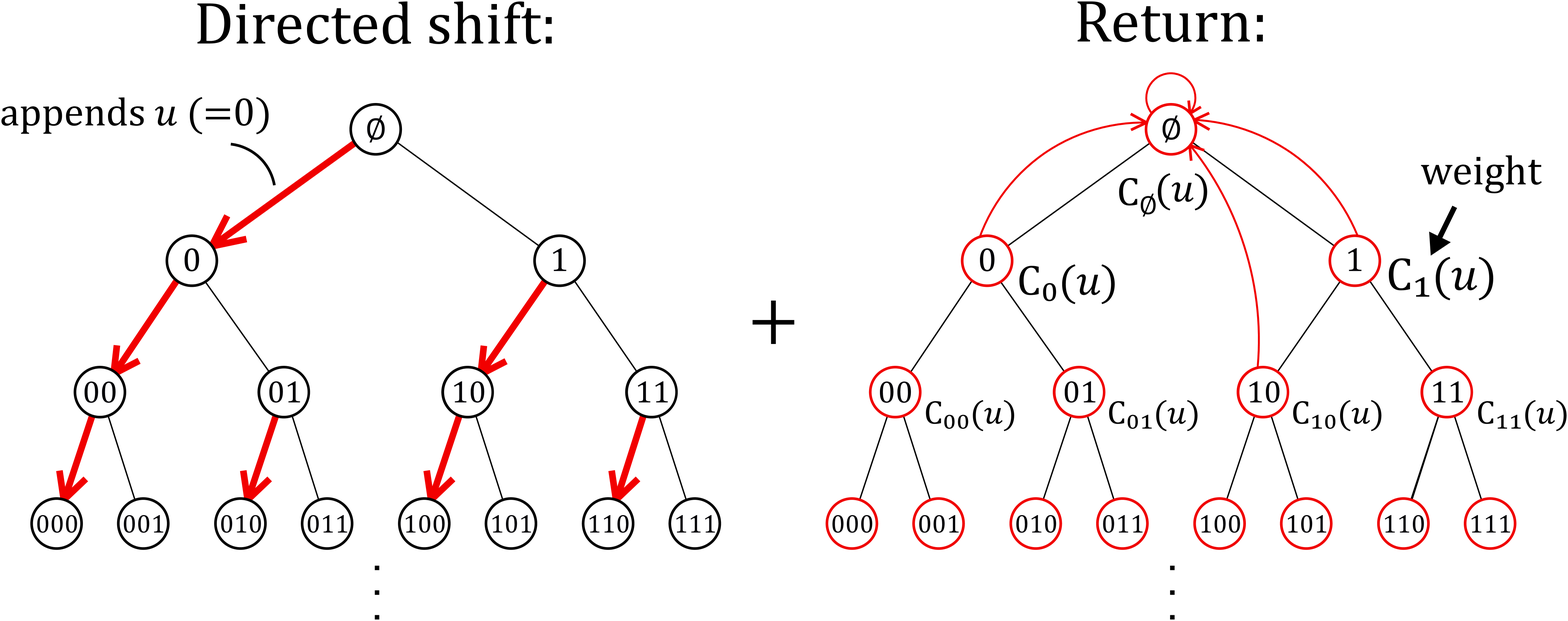}
\vspace{3pt}%
\caption{A description of the operator $M_u$ when $I_1 =\{0, 1\}$ and $u = 0 \in J$} \label{figure: operator in 3d}
\end{figure}

Let $\Gamma$ be the set of all finite words formed from the letters in $J$, including the null word, which we interpret as a rooted tree:
\[ \Gamma = \{\varnothing\} \cup \bigcup_{k = 1}^\infty J^k. \]
Let $\mathbb{R}^{\bigoplus \Gamma}$ denote the direct sum of $\mathbb{R}$ indexed by $\Gamma$, where an element has only finitely many non-zero entries. For $u \in J$, define a linear operator $M_u: \mathbb{R}^{\bigoplus \Gamma} \rightarrow \mathbb{R}^{\bigoplus \Gamma}$ as the sum of a directed shift and a return map (see Figure \ref{figure: operator in 3d}):
\begin{equation*}
\left( M_u \left( (x_{\mu})_{\mu \in \Gamma} \right) \right)_\lambda
=
\begin{dcases}
\sum_{s \in \Gamma} C_s(u) x_s \text{\quad ( if $\lambda = \varnothing$. )} \\
x_{\lambda'}^{} \text{\quad ( if $\lambda = \lambda' \spa u$, the concatenation of $\lambda'$ and $u$. )} \\
0 \text{\quad ( otherwise. )}
\end{dcases}
\end{equation*}
When $u$ is removable ($u \notin J$), $M_u$ is defined similarly, although without the shift operator as
\begin{equation*}
\left( M_u \left( (x_{\mu})_{\mu \in \Gamma} \right) \right)_\lambda
=
\begin{dcases}
\sum_{s \in \Gamma} C_s(u) x_s \text{\quad ( if $\lambda = \varnothing$. )} \\
0 \text{\quad ( otherwise. )}
\end{dcases}
\end{equation*}

Define $\Phi_0 \in \mathbb{R}^{\Gamma}$ by $(\Phi_0)_{\varnothing} = 1$ and $0$ elsewhere. Let $A = \sum_{w \in I_1 \times I_2} A_w$.

\begin{theorem} \label{theorem: tower decomposition for 3d}
Suppose that $X$ has a recursive structure with $\boldsymbol{v} \in \mathbb{R}^n$, and that $A$ is primitive. Then,
\begin{equation*}
\mathrm{dim}_{\mathrm{H}}(X)
=
\lim_{N \to \infty} \spa \frac{1}{N} \log_{m_1^{}}{\hspace{-3pt}\sum_{(u_1^{}, \ldots, u_N^{}) \in I_1^N}
{\norm{ M_{u_N^{}} \spa M_{u_{N-1}^{}} \cdots \spa M_{u_1^{}} \spa \Phi_0 }_1
}^{a_1^{}}}.
\end{equation*}
\end{theorem}

\begin{proof}

The proof proceeds by constructing a tower decomposition of the summand over the tree $\Gamma$ and expressing it as a composition of linear operators $M_u$.

Let $N$ be a natural number and $s = (s_1^{}, \ldots, s_N^{}) \in I_1^N$. We begin by defining $\Phi_N(s) \in \mathbb{R}^{\Gamma}$. For a natural number $k$, define $s[k]$ to be the last $k$ components of $s$ from the end, that is, $s[k] = (s_{N-k+1}^{}, \ldots, s_{N}^{})$. Also, $s|_k$ is the first $k$ components; $s|_k = (s_1^{}, \ldots, s_k^{})$. Set $s[0] = \varnothing$. Let
\[ \left( \Phi_N(s) \right)_{s} = \left[ \boldsymbol{v}^{\top} \spa \mathbb{e} \right]^{a_2^{}}, \]
\[ \left( \Phi_N(s) \right)_{s[N-1]}
=
\sum_{p \in Q(s_{1}^{})} \hspace{2.5pt} \left[ \boldsymbol{v}^{\top} \spa A_p \spa \mathbb{e} \right]^{a_2^{}}, \]
and for $0 \leq k \leq N-2$,
\[ \left( \Phi_N(s) \right)_{s[k]}
=
\sum_{p \in Q(s_{N-k}^{})} \hspace{2.5pt}
\sum_{(w_i^{})_i \in W^{N-k-1}(s|_{N-k-1})} \left[ \boldsymbol{v}^{\top} \spa A_{w_1^{}} \cdots A_{w_{N-k-1}^{}} \spa A_{p} \spa \mathbb{e} \right]^{a_2^{}}. \]
Here, $\mathbb{e}$ denotes the column vector with $1$ in every entry. We define all other entries of $\Phi_N(s)$ to be $0$. $\big($ Although the above definition depends only on $k$ and not on $s[k]$, we will need this distinction for later discussion.$\big)$

\begin{claim} \label{claim: tower decomposition for 3d}
We have the following relation:
\begin{equation*}
\mathrm{dim}_{\mathrm{H}}(X)
=
\lim_{N \to \infty} \frac{1}{N} \log_{m_1^{}}
\sum_{s \in I_1^N} \norm{ \Phi_N(s) }_1^{a_1}.
\end{equation*}
\end{claim}

\begin{proof}[Proof of Claim \ref{claim: tower decomposition for 3d}]
Since $A = \sum_{w \in I_1 \times I_2} A_w \in \mathrm{M}_n(\mathbb{Z}_{\geq0})$ is primitive, there is an integer $d$ such that every entry of $A^d$ is at least 1. For a natural number $k$, take a string $w^{(k)} = (w_1^{(k)}, \ldots, w_k^{(k)}) \in I_2^k$ such that $A_{w_1^{(k)}} \cdots A_{w_k^{(k)}} \ne O$. Then,
\[ \norm{\boldsymbol{x}^{\top} }_1 \leq \norm{ \boldsymbol{x}^{\top} \spa A^{d} \spa A_{w_1^{(k)}} \cdots A_{w_k^{(k)}} }_1 \]
for every $\boldsymbol{x} \in (\mathbb{R}_{\geq 0})^n$. Using the fact that $(x + y)^{\alpha} \leq x^{\alpha} + y^{\alpha}$ for $x, y \geq 0$ when $0 \leq \alpha \leq 1$, we can evaluate as follows.
\begin{flalign*}
& \sum_{s \in I_1^N} \norm{ \Phi_N(s) }_1^{a_1}&
\end{flalign*}
\begin{flalign*} 
&\leq
\sum_{s \in I_1^N} \hspace{-2pt}
\left( \hspace{-2pt}
(\Phi_N(s))_s +
\sum_{k = 0}^{N-1} \hspace{0.5pt}
\sum_{p \in Q(s_{N-k}^{})} \hspace{1pt}
\sum_{\substack{
	(w_i)_i \spa \in \\
	W_2^{N-k-1}(s|_{N-k-1})
	}} \hspace{-1pt}
\left[ \boldsymbol{v}^{\top} \spa A_{w_1^{}} \cdots A_{w_{N-k-1}^{}} \spa A_{p} \spa A^d \spa A_{w_1^{(k)}} \cdots A_{w_k^{(k)}} \spa \mathbb{e} \right]^{a_2}
\right)^{\hspace{-2pt} a_1^{}} &
\end{flalign*}
\begin{flalign*} 
&=
\sum_{s \in I_1^N}
\left(
\left[ \boldsymbol{v}^{\top} \spa \mathbb{e} \right]^{a_2^{}} +
\sum_{k} \hspace{2pt}
\sum_{p} \hspace{2.5pt}
\sum_{(w_i)_i}
\left[ \boldsymbol{v}^{\top} \spa A_{w_1^{}} \cdots A_{w_{N-k-1}^{}} \spa A_{p} \spa \left( \sum_{w \in I_1 \times I_2} A_w \right)^d \spa A_{w_1^{(k)}} \cdots A_{w_k^{(k)}} \spa \mathbb{e} \right]^{a_2} \hspace{-1pt}
\right)^{a_1^{}} &
\end{flalign*}
\begin{flalign*} 
& \leq
{(m_1m_2)}^{d+1} \spa N \spa \sum_{s \in I_1^{N+d}}
\left(
\sum_{(w_i)_i \in W_2^{N+d}(s)}
{\left[ \boldsymbol{v}^{\top} \spa A_{w_1^{}} \cdots A_{w_{N+d}^{}} \spa \mathbb{e}
\right]
}^{a_2^{}}
\right)^{a_1^{}} &
\end{flalign*}
\begin{flalign*} 
& \leq
{(m_1m_2)}^{d+1} \spa N \spa \norm{\boldsymbol{v}}_{\infty}
\sum_{s \in I_1^{N+d}}
\left(
\sum_{(w_i)_i \in W_2^{N+d}(s)}
{
\left[ \mathbb{e}^{\top} \spa A_{w_1^{}} \cdots A_{w_{N+d}^{}} \spa \mathbb{e}
\right]
}^{a_2^{}}
\right)^{a_1^{}} &
\end{flalign*}
\begin{flalign*} 
& \leq
{(m_1m_2)}^{d+1} \spa N \spa \norm{\boldsymbol{v}}_{\infty}
\sum_{s \in I_1^{N+d}}
\left(
\sum_{p \in Q(0)}
\sum_{(w_i)_i \in W_2^{N+d}(s)}
{
\left[ \boldsymbol{v}^{\top} \spa A^d \spa A_{w_1^{}} \cdots A_{w_{N+d}^{}} \spa A^d \spa A_p \spa \mathbb{e}
\right]
}^{a_2^{}}
\right)^{a_1^{}} &
\end{flalign*}
\begin{flalign*} 
& \leq
{(m_1m_2)}^{3d+1} \spa N \spa \norm{\boldsymbol{v}}_{\infty}
\sum_{s \in I_1^{N+3d+1}} \norm{ \Phi_N(s) }_1^{a_1}. &
\end{flalign*}
Taking the logarithm and letting $N \rightarrow \infty$ finishes the proof.
\end{proof}

Now, take $s \in I_1^N$ and $u \in I_1$, and consider the concatenation $su \in I_1^{N+1}$. We see from the definition that for a non-negative integer $k$,
\[ \left( \Phi_{N+1}(su) \right)_{s[k] \spa u} = \left( \Phi_{N+1}(su) \right)_{(su)[k+1]} = \left( \Phi_{N}(s) \right)_{s[k]}. \]
Also, the vector $\boldsymbol{v}^{\top} \spa A_{w_1^{}} \cdots A_{w_{N-k-1}^{}} \spa A_p$ in the definition of $\Phi_N(s)$ is already a constant multiple of $\boldsymbol{v}^{\top}$. Therefore, we have for any $0 \leq k \leq N$ and $u \in I_1$,
\begin{flalign*}
&\; & \sum_{q \in Q(u)} \hspace{3.5pt} \sum_{(z_i^{})_i \in R^k(s[k])} \hspace{3.5pt} \sum_{p \in Q(s_{N-k}^{})} \hspace{3.5pt} \sum_{(w_i^{})_i \in W_2^{N-k-1}(s|_{N-k-1})}
\left[ \boldsymbol{v}^{\top} A_{w_1^{}} \cdots A_{w_{N-k-1}^{}} \spa A_{p} \spa A_{z_1^{}} \cdots A_{z_k^{}} \spa A_{q} \spa \mathbb{e} \right]^{a_2^{}} &
\end{flalign*}
\begin{flalign*}
& \phantom{\hspace{25pt}} = \sum_{q \in Q(u)} \hspace{3.5pt} \sum_{w \in R^k(s[k])} {D_w(q)}^{a_2} \spa \left( \Phi_N(s) \right)_{s[k]} & \\
& \phantom{\hspace{25pt}} = C_{s[k]}(u) \left( \Phi_N(s) \right)_{s[k]}. &
\end{flalign*}
Since $W_2^N(s) = \bigcup_{k = 0}^N W_2^{N-k-1}(s|_{N-k-1}) \times Q(s_{N-k}^{}) \times R^k (s[k]) \footnote{Here, the sets $W_2^{-1}$, $Q(s_0^{})$, and $R^0$ are regarded as empty.}$ is a disjoint partition (by considering the last index, say $N-k$, that multiplies $A_p$ with some $p \in Q(s_{N-k}^{})$,) we conclude
\[ \left( \Phi_{N+1}(su) \right)_{\varnothing} = \sum_{\lambda \in \Gamma} C_{\lambda}(u) \spa \left( \Phi_N(s) \right)_\lambda. \]
Here, we used the fact that $\left( \Phi_N(s) \right)_\lambda = 0$ for all $\lambda$ except for $\lambda = s[k]$ with some $k$. By the definition of $M_u$, we have
\[ \Phi_{N+1}(su) = M_u \spa \Phi_N(s). \]
We conclude that for any string $s = (s_1^{}, \ldots, s_N^{}) \in I_1^N$,
\[ \Phi_N(s) = M_{s_N^{}} \cdots M_{s_1^{}} \spa \Phi_0. \]

We combine this with Claim \ref{claim: tower decomposition for 3d} and complete the proof.

\end{proof}

When there is a removable index $u \in I_1$, the operator $M_u$ has a $1$-dimensional image (since it does not include the ``shift'' operator), leading to an intriguing result similar to the one in Theorem \ref{theorem: tower decomposition for 2d}.

\begin{proposition} \label{proposition: reapplication}
Suppose that a sofic set $X \subset \mathbb{T}^3$ satisfies the following conditions. \\
(1) $X$ has a recursive structure with $\boldsymbol{v} \in \mathbb{R}^n$, and $A$ is primitive. \\
(2) There is at least one non-removable index in $I_1$\footnote{Otherwise, the dimension is easily determined.}. \\
(3) There is a removable index $u \in I_1$, and a string $(t_1^{}, \ldots, t_L^{}) \in I_1^L$ such that $M_u \spa M_{t_L^{}} \cdots M_{t_1^{}}$ is increasing with respect to the $l^1$-norm\footnote{That is, $\norm{M_u \spa M_{t_L^{}} \cdots M_{t_1^{}} \spa \boldsymbol{x}}_1 \geq \norm{\boldsymbol{x}}_1$ for any $\boldsymbol{x} \in \mathbb{R}^{\bigoplus \Gamma}$. This is easily satisfied in most cases.}. \\
Let $s = (s_1^{}, \ldots, s_k^{}) \in I_1^k$ be a string, and define $c(s) \geq 0$ by
\[ M_u \spa M_{s_{N}^{}} \cdots M_{s_1^{}} \spa \Phi_0 = c(s) \spa \Phi_0. \]
Next, define $b_k$ by
\[ b_k = \sum_{s \in (I_1 \backslash \{u\})^k} c(s)^{a_1}. \]
We then conclude that
\[ \mathrm{dim}_{\mathrm{H}}(X) = \log_{m_1^{}} r, \]
where $r$ is the unique positive solution to the equation
\[ r = b_0 + \frac{b_1}{r} + \frac{b_2}{r^2} + \cdots. \]
\end{proposition}

\begin{proof}
Without loss of generality, assume that $u = 0$. We begin by constructing the tower decomposition $\Psi_N = \left( \Psi_N(k) \right)_{k = 0}^\infty \in \mathbb{R}^{\mathbb{N}_0}$ for a natural number $N$. First, set $\Psi_0 = (1, 0, 0, \ldots)^{\top}$, and let for $0 \leq k \leq N-1$
\[ \Psi_N(k) =
\sum_{(s_1^{}, \ldots, s_{N-k-1}^{}) \in I_1^{N-k-1}}
\norm{ M_0 \spa M_{s_{N-k-1}^{}} \cdots M_{s_1^{}} \spa \Phi_0 }_1^{a_1^{}}. \]
The following evaluation follows from assumption (2) and (3).
\begin{align*}
\norm{ \Psi_N }_1
\leq 
& \sum_{k = 0}^{N-1} \hspace{2pt} \sum_{(w_1^{}, \ldots, w_k^{}) \in I_1^{k}} \hspace{2pt} \sum_{(s_1^{}, \ldots, s_{N-k-1}^{}) \in I_1^{N-k-1}}
\norm{ M_{w_k^{}} \cdots M_{w_1^{}} \spa M_0 \spa M_{t_L^{}} \cdots M_{t_1^{}} \spa M_0 \spa M_{s_{N-k-1}^{}} \cdots M_{s_1^{}} \spa \Phi_0 }_1^{a_1^{}}. & \\
\leq
& \spa N \sum_{(s_1^{}, \ldots, s_{N+L+2}^{}) \in I_1^{N+L+2}}
\norm{ M_{s_{N+L+2}^{}} \cdots M_{s_1^{}} \spa \Phi_0 }_1^{a_1^{}} \\
\leq
& \spa N \sum_{(s_1^{}, \ldots, s_{N+L+2}^{}) \in I_1^{N+L+2}}
\norm{ M_0 \spa M_{t_L^{}} \cdots M_{t_1^{}} \spa M_{s_{N+L+2}^{}} \cdots M_{s_1^{}} \spa \Phi_0 }_1^{a_1^{}} \\
\leq
& \spa N \spa \norm{\Psi_{N+2L+2}}_1.
\end{align*}
The estimation above implies

\[
\lim_{N \to \infty} \spa \frac{1}{N} \log_{m_1^{}}{\hspace{-3pt}\sum_{(s_1^{}, \ldots, s_N^{}) \in I_1^N}
{\norm{ M_{s_N^{}} \spa M_{s_{N-1}^{}} \cdots \spa M_{s_1^{}} \spa \Phi_0 }_1
}^{a_1^{}}}
=
\lim_{N \to \infty} \spa \frac{1}{N} \log_{m_1^{}} \norm{\Psi_N}_1. \]

Next, we define a linear operator $L: \mathbb{R}^{\bigoplus \mathbb{N}_0} \rightarrow \mathbb{R}^{\bigoplus \mathbb{N}_0}$ by
\begin{equation*}
L = 
\begin{pmatrix}
b_0 & b_1 & b_2 & \cdots \\
1 & 0 & 0 & \\
0 & 1 & 0 & \cdots \\
0 & 0 & 1 & \\
& \vdots & & \ddots \\
\end{pmatrix}.
\end{equation*}
Then, $\Psi_{N+1} = L \spa \Psi_N$, and $\Psi_N = L^N \Psi_0$. The remainder of the proof follows in the same way as in Theorem \ref{theorem: tower decomposition for 2d}.
\end{proof}

\section{Examples} \label{section: examples}

In this section, we provide a detailed explanation to the examples mentioned in the introduction, starting with the planar cases.

\begin{example}

Let $I = \{0, 1\} \times \{0, 1, 2\}$. Consider the directed graph in Figure \ref{figure: digraph for planar 01 redrawn} labeled with $I$.
\begin{figure}[h!]
\includegraphics[width=\textwidth-5.5cm]{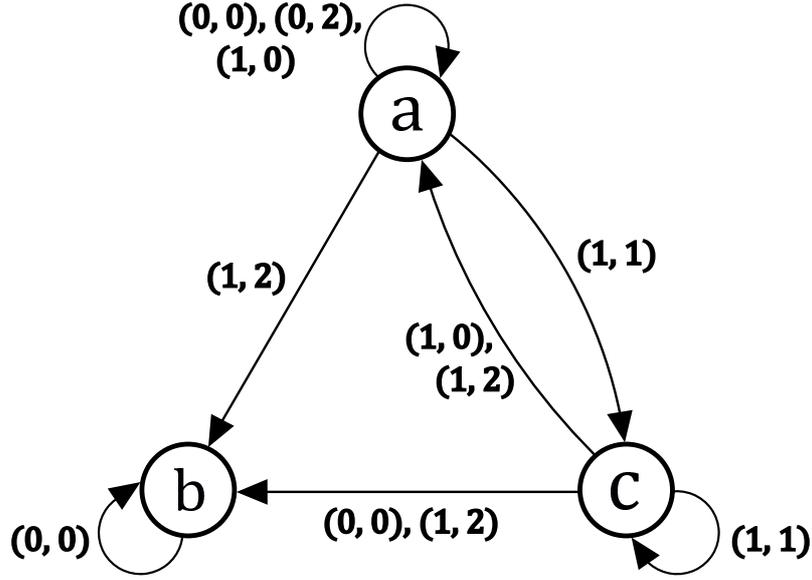}
\caption{A digraph labeled with $I = \{0, 1\} \times \{0, 1, 2\}$} \label{figure: digraph for planar 01 redrawn}
\end{figure}
Then, the adjacency matrices are
\begin{equation*}
A_0 =
\begin{pmatrix}
2 & 0 & 0\\
0 & 1 & 0\\
0 & 1 & 0
\end{pmatrix},
\hspace{4pt} A_1 =
\begin{pmatrix}
1 & 1 & 1 \\
0 & 0 & 0 \\
2 & 1 & 1
\end{pmatrix}.
\end{equation*}
The summed matrix $A = A_0 + A_1$ is primitive. Furthermore,
\begin{equation*}
A_0 \spa A_1 =
\begin{pmatrix}
2 & 2 & 2\\
0 & 0 & 0\\
0 & 0 & 0
\end{pmatrix}
\end{equation*}
has a $1$-dimensional image, $\mathrm{Span}\{(1, 1, 1)\}$, satisfying the condition in Theorem \ref{theorem: tower decomposition for 2d}. Excluding the string ``$01$'', the only possible strings are of the form $11\cdots100\cdots0$. For such a string $(u_1^{}, \ldots, u_N^{})$, we have for some $k$
\begin{flalign*}
(1, 1, 1) \spa A_{u_1^{}} \cdots A_{u_N^{}} \spa A_0 \spa A_1
& = (1, 1, 1) \spa A_1^{k} \spa A_0^{N-k} \spa A_0 \spa A_1 & \\
& = 2^{N-k+1} \Big( ( 2 + 2\sqrt{2} )(1+ \sqrt{2})^k - ( 2 - 2\sqrt{2} ) (1- \sqrt{2})^k \Big) \spa (1, 1, 1). &
\end{flalign*}
Let
\begin{align*}
C_{N, k} = 2^{N-k+1} \Big( ( 2 + 2\sqrt{2} )(1+ \sqrt{2})^k - ( 2 - 2\sqrt{2} ) (1- \sqrt{2})^k \Big).
\end{align*}
Then, by Theorem \ref{theorem: tower decomposition for 2d}, $\mathrm{dim}_{\mathrm{H}}(X) = \log_2{r} = 1.6416\cdots$, where $r = 3.1201\cdots$ satisfies
\begin{equation*}
r^2 = \sum_{N =1}^\infty \left( \sum_{k = 0}^N
{C_{N, k}}^{\log_3{2}} \right) r^{-N}.
\end{equation*}

\end{example}

\spa

\begin{example}

Let $G$ be a directed graph with $2$ vertices, and $I = \{0, 1, 2\} \times \{0, 1, 2, 3\}$. Consider a sofic system with the following adjacency matrices.
\begin{equation*}
A_0 =
\begin{pmatrix}
1 & 0 \\
2 & 0 \\
\end{pmatrix},
\hspace{4pt} A_1 =
\begin{pmatrix}
2 & 1 \\
1 & 2 \\
\end{pmatrix},
\hspace{4pt} A_2 =
\begin{pmatrix}
3 & 2 \\
2 & 3 \\
\end{pmatrix}.
\end{equation*}
In this case, $A_0$ has a $1$-dimensional image, $\mathrm{Span}\{ (1, 0) \}$. Furthermore, $A_1$ and $A_2$ commute, which implies that for every string $(u_1, \ldots, u_N) \in \{1, 2\}^N$,
\begin{align*}
A_{u_1^{}} \cdots A_{u_N^{}} = A_1^k \spa A_2^{N-k}
\end{align*}
with some $k$. Then, \vspace{-8pt}
\begin{align*}
(1, 0) \spa A_{u_1^{}} \cdots A_{u_N^{}} \spa A_0
&= \frac{3^k \spa 5^{N-k} -1}{2} (1, 0).
\end{align*}
Counting the number of strings with such a $k$, we have
\begin{equation*}
C_k = \sum_{k = 0}^N
\begin{pmatrix}
N \\
k
\end{pmatrix}
\left( \frac{3^k \spa 5^{N-k} -1}{2} \right)^{\log_4{3}}.
\end{equation*}
Then, $\mathrm{dim}_{\mathrm{H}}(X) = \log_3{r} = 1.6994\cdots$, where $r$ satisfies
\begin{align*}
r = \sum_{N=0}^{\infty} \left(\sum_{k = 0}^N
\begin{pmatrix}
N \\
k
\end{pmatrix}
\left( \frac{3^k \spa 5^{N-k} -1}{2} \right)^{\log_4{3}} \right) r^{-N} = 6.4693\cdots.
\end{align*}

\end{example}

Next, we provide a detailed explanation of the examples of sofic sets in $\mathbb{T}^3$ introduced earlier.

\begin{example}

Let $I = \{0, 1\} \times \{0, 1, 2\} \times \{0, 1, 2, 3\}$. Consider the directed graph in Figure \ref{figure: digraph for 3d example reappear} labeled with $I$.
\begin{figure}[h!]
\includegraphics[width=\textwidth-4cm]{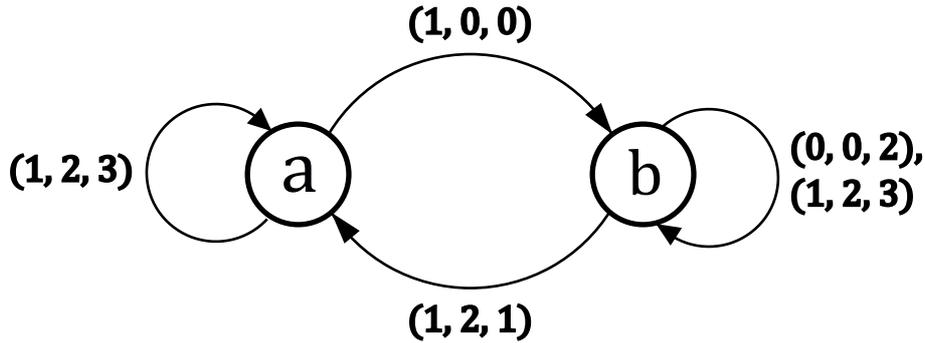}
\caption{A digraph labeled with $I = \{0, 1\} \times \{0, 1, 2\} \times \{0, 1, 2, 3\}$} \label{figure: digraph for 3d example reappear}
\end{figure} \\
Then, the adjacency matrices are
\begin{equation*}
\hspace{4pt}
A_{(0, 0)} =
\begin{pmatrix}
0 & 0 \\
0 & 1 \\
\end{pmatrix},
\hspace{4pt} A_{(1, 0)} =
\begin{pmatrix}
0 & 1 \\
0 & 0 \\
\end{pmatrix},
\hspace{4pt} A_{(1, 2)} =
\begin{pmatrix}
1 & 0 \\
1 & 1 \\
\end{pmatrix}.
\end{equation*}
This sofic set has a recursive structure with $(0, 1)$. Moreover, $0 \in I_1$ is removable, so we have $J = \{1\}$. The coefficients are; $C_\varnothing (0) = 1$, $C_s(0) = 1$ for $s \in \{ 1 \}^N$, $C_\varnothing (1) = 0$, and $C_s(1) = N^{\log_4{3}}$ for $s \in \{ 1 \}^N$. The operator $M_1$ acts on $\mathbb{R}^{\bigoplus \mathbb{N}_0}$, and it is given by
\begin{equation*}
M_1 =
\begin{pmatrix}
0 & 1 & 2^{\log_4{3}} & 3^{\log_4{3}} & 4^{\log_4{3}} & \cdots \\
1 & 0 & 0 & \cdots \\
0 & 1 & 0 & 0 & \cdots \\
0 & 0 & 1 & 0 & 0 & \cdots \\
& \vdots & & & \ddots \\
\end{pmatrix}.
\end{equation*}
Let $\Psi_0 = (1, 0, 0, \ldots)^{\top}$. We have for any natural number $k$
\begin{equation*}
b_k =
\norm{
M_1^k \Phi_0
}_{1}^{\log_3{2}}.
\end{equation*}
Then, $\mathrm{dim}_{\mathrm{H}}(X) = \log_2{r} = 1.1950\cdots$, where $r$ satisfies
\begin{align*}
r
= \sum_{k = 0}^{\infty} \frac{b_k}{r^k}
&= 1 + \frac{1}{r} + \frac{\sqrt{2}}{r^2} + \frac{ (2 + 2^{\log_4{3}})^{\log_3{2}} }{r^3} + \frac{ (3 + 2^{\log_4{3}} + 3^{\log_4{3}})^{\log_3{2}} }{r^4} + \cdots & \\
&  = 2.2894\cdots. &
\end{align*}

\end{example}

\begin{example}

Let $I = \{0, 1, 2\} \times \{0, 1, 2, 3\} \times \{0, 1, 2, 3, 4\}$. Consider the directed graph in Figure \ref{figure: digraph for 3d example 2 redrawn} labeled with $I$.
\begin{figure}[h!]
\includegraphics[width=\textwidth-4cm]{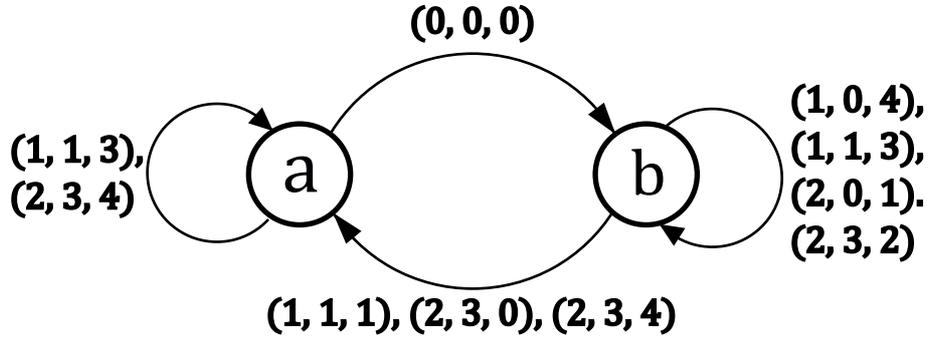}
\vspace{3pt}%
\caption{A digraph labeled with $I = \{0, 1, 2\} \times \{0, 1, 2, 3\} \times \{0, 1, 2, 3, 4\}$} \label{figure: digraph for 3d example 2 redrawn}
\end{figure} \\
Then, the adjacency matrices are
\begin{equation*}
A_{(0, 0)} =
\begin{pmatrix}
0 & 1 \\
0 & 0 \\
\end{pmatrix},
\hspace{4pt} A_{(1, 0)} =
\begin{pmatrix}
0 & 0 \\
0 & 1 \\
\end{pmatrix},
\hspace{4pt} A_{(1, 1)} =
\begin{pmatrix}
1 & 0 \\
1 & 1 \\
\end{pmatrix},
\hspace{4pt} A_{(2, 0)} =
\begin{pmatrix}
0 & 0 \\
0 & 1 \\
\end{pmatrix},
\hspace{4pt} A_{(2, 3)} =
\begin{pmatrix}
1 & 0 \\
2 & 1 \\
\end{pmatrix}.
\end{equation*}

Let $a_1^{} = \log_4{3}$ and $a_2^{} = \log_5{4}$. This sofic set has a recursive structure with $(0, 1)$, and $J = \{1, 2\}$ since $0 \in I_1$ is removable. The coefficients for $u = 0$ are; $C_\varnothing (0) = 0$, $C_s(0) = (N + k)^{a_2^{}}$ for $s \in J^N$ that has $k$ number of $2$s. Also, $C_s(1) = C_s(2) = 1$ for every $s \in \Gamma$. Then, by considering
\[ M_u \spa M_{s_{N}^{}} \cdots M_{s_1^{}} \spa \Phi_0, \]
we see that
\begin{align*}
\hspace{10pt} b_N
&= \sum_{(s_1, \ldots, s_N) \in \{1, 2\}^N} \Bigg( \big( N + \#\left\{ j \setcond 1 \leq j \leq N, s_j^{} =2\right\} \big)^{a_2^{}} \hspace{20pt} &\\
&\hspace{150pt} + \sum_{k = 1}^{N-1} 2^{k-1} \big( N + \#\left\{j \setcond k \leq j \leq N \text{\hspace{1pt} and \hspace{1pt}} s_j^{} = 2\right\} \big)^{a_2^{}} \Bigg)^{a_1^{}}. &
\end{align*}
We conclude that $\mathrm{dim}_{\mathrm{H}}(X) = \log_2{r} = 2.224\cdots$, where $r$ satisfies
\begin{align*}
r
= \sum_{k = 0}^{\infty} \frac{b_k}{r^k}
= 4.673\cdots.
\end{align*}

\end{example}

\section*{Acknowledgement}
I am deeply grateful to my mentor, Masaki Tsukamoto, whose expertise and insightful feedback have been instrumental in shaping both this paper and my understanding of ergodic theory. 

I also want to acknowledge my family and friends for their unwavering support, which has been a foundation throughout my life.

This paper owes much to the collective support and intellectual environment of the academic community I have been fortunate to be a part of.

\vspace{0.5cm}

\address{
Department of Mathematics, Kyoto University, Kyoto 606-8501, Japan}

\textit{E-mail}: \texttt{alibabaei.nima.28c@st.kyoto-u.ac.jp}

\end{document}